\documentclass[10pt]{article}
\usepackage{a4wide}
\usepackage{latexsym}
\usepackage{amsmath}
\usepackage{amsfonts}
\usepackage{amssymb}
\usepackage{amsthm}
\usepackage[all]{xy}

\def\cH{\mathcal{H}}
\def\cL{\mathcal{L}}
\def\cR{\mathcal{R}}

\def\cK{\mathcal{K}}
\def\cF{\mathcal{F}}

\def\grp{^\#}

\def\iff{\Leftrightarrow}

\def\beqn{\begin{eqnarray*}}
\def\eeqn{\end{eqnarray*}}
\def\barr{\begin{array}}
\def\earr{\end{array}}

\def\ben{\begin{enumerate}}
\def\een{\end{enumerate}}
\def\l({\left(}
\def\r){\right)}

\def\bmx{\left[\begin{array}}
\def\emx{\end{array}\right]}

\newtheorem{theorem}{Theorem}[section]
\newtheorem{lemma}[theorem]{Lemma}
\newtheorem{proposition}[theorem]{Proposition}
\newtheorem{corollary}[theorem]{Corollary}
\newtheorem{definition}[theorem]{Definition}

\begin{document}

\title{Classes of semigroups modulo Green's relation $\cH$}

\author{Xavier Mary\footnote{email: xavier.mary@u-paris10.fr}\\
\textit{\small Universit\'e Paris-Ouest Nanterre-La D\'efense, Laboratoire Modal'X}}

\date{}

\maketitle


\begin{abstract}
Inverses semigroups and orthodox semigroups are either defined in terms of inverses, or in terms of the set of idempotents $E(S)$. In this article, we study analogs of these semigroups 
defined in terms of inverses modulo Green's relation $\cH$, or in terms of the set of group invertible elements $\cH(S)$, that allows a study of non-regular semigroups.  We then study the interplays between these new classes of semigroups, as well as with known classes of semigroups (notably inverse, orthodox and cryptic semigroups).\end{abstract}

\section{Introduction}


The study of special classes of semigroups relies in many cases on properties of the set of idempotents, or of  regular pairs of elements. In \cite{mary} a special weak inverse, called inverse along an element, was introduced, and in \cite{maryPatricio}, it was interpreted as a kind of inverse modulo $\cH$. The aim of this article is to study special classes of semigroups defined in terms of these inverses modulo Green's relation $\cH$, or in terms of the set of group invertible elements, that can be seen as idempotents modulo $\cH$. Relations with known classes of semigroups will be investigated.\\

In this paper, $S$ is a semigroup and $S^1$ denotes the monoid generated by $S$. $E(S)$ denotes the set of idempotents, and $Z(E(S))=\{x\in S,\, xe=ex\; \forall e\in E(S)\}$ its centralizer.

We say $a$ is (von Neumann) regular in $S$ if $a\in aSa$. A particular solution to $axa=a$ is called an associate, or inner inverse, of $a$.  A solution to $xax=a$ is called a weak (or outer) inverse. Finally, an element that satisfies $axa=a$ and $xax=x$ is called an inverse (or reflexive inverse, or relative inverse) of $a$. The set of all associates of $a$ is denoted by $A(a)$ and the set of all inverses of $a$ by $V(a)$. If $a'\in V(a)$, we also say that $(a, a')$ is a regular pair.

A commuting inverse, if it exists, is unique and denoted by $a\grp$. It is usually called the group inverse of $a$. We let $\cH(S)$ denote the set of group invertible elements. \\

We will make use of the Green's preorders and relations in a semigroup \cite{Green51}. For elements $a$ and $b$ of $S$, Green's preorders  $\leq_{\cL}$, $\leq_{\cR}$ and
$\leq_{\cH}$ are defined by
\begin{align*}
a \leq_{\cL} b&\Longleftrightarrow S^1 a\subset S^1 b\Longleftrightarrow \exists x\in S^1,\;  a = xb;\\ 
a \leq_{\cR} b& \Longleftrightarrow  aS^1\subset  bS^1\Longleftrightarrow \exists x\in S^1,\;  a = bx;\\ 
a \leq_{\cH} b&\Longleftrightarrow \{a\leq_{\cL} b\text{ and }a\leq_{\cR} b\}.
\end{align*}

If $\leq_{\cK}$ is one of these preorders, then $a\cK b\iff \{a\leq_{\cK}b \text{ and } b\leq_{\cK}a\}$, and $\cK_a=\{b\in S,\; b\cK a\}$ denotes the $\cK$-class of $a$. 

We will use the following classical lemmas. Let $a,b,c\in S$.

\begin{lemma}[Cancellation]\label{lemma-cancel}
$\,$\\
\vspace{-.5cm}
\begin{align*}
a \leq_{\cL} b &\Rightarrow \{\forall x,y\in S^1,\; bx=by \Rightarrow ax=ay\};\\ 
a \leq_{\cR} b &\Rightarrow \{\forall x,y\in S^1,\; xb=yb \Rightarrow xa=ya \}. 
\end{align*}
\end{lemma}

\begin{lemma}\label{lemma-elementary}
$\,$\\
\vspace{-1cm}
\begin{align*}
ca \leq_{\cL} a &,\; ac \leq_{\cR} a,\; aca\leq_{\cH} a;\\ 
a \leq_{\cL} b &\Rightarrow ac \leq_{\cL} bc \text{ (Right congruence)};\\ 
a \leq_{\cR} b &\Rightarrow ca \leq_{\cR} cb \text{ (Left congruence)}. 
\end{align*}
\end{lemma}

Note that $\cH$ is not a congruence in general.

\begin{lemma}\label{lemma-L}
$\,$\\
\vspace{-.5cm}
\ben
\item Let $b$ be a regular element. Then $a\leq_{\cL} b\iff a=ab'b$ for one (any) inverse $b'\in V(b)$;
\item Assume $S$ is inverse. Then $a\cL b\iff a^{-1}a=b^{-1}b$.
\een
\end{lemma}

The purpose of this article is to study inverses modulo Green's relations $\cH$.

\begin{definition}
We call a particular solution to $axa\cH a$ an associate of $a$ modulo $\cH$.  If also $xax\cH a$, then $x$ is called an inverse of $a$ modulo $\cH$. Finally, we denote the set of all associates of $a$ modulo $\cH$ by $A(a)[\cH]$ , and the set of inverses of $a$ modulo $\cH$ by $V(a)[\cH]$.
\end{definition}

We recall the following characterization of group invertibility in terms of Green's relation $\cH$ (see \cite{Green51}, \cite{Clifford56}) and inverses:
\begin{lemma}\label{lemma_Clifford}
$\,$\\
\vspace{-.5cm}
\ben
\item $a\grp$ exists if and only if $a\cH a^2$ if and only if $\cH_a$ is a group.
\item Let $(a,a')$ be a regular pair. Then $aa'=a'a$ if and only if $a\cH a'$. 
\een
\end{lemma}

From this lemma, we get that the set of idempotents modulo $\cH$ is precisely the set of group invertible (completely regular) elements, and can be characterized as the union of all (maximal) subgroups of $S$: $$E(S)[\cH]=\cH(S)=\bigcup_{e\in E(S)}\cH_e.$$

\section{Invertibility modulo $\cH$}

In \cite{mary} a special weak inverse, called inverse along an element, was introduced, and in \cite{maryPatricio}, it was interpreted as a kind of inverse modulo $\cH$. We recall the definition and properties of this inverse, and refer to \cite{mary} and \cite{maryPatricio} for the proofs.

\begin{definition}
Given $a,a'$ in $S$, we say $a'$ is invertible along $a$ if there exists $b\in S$ such that  $ba'a=a=aa'b$ and $b\leq_\cH a$. If such an element exists then it is unique and is denoted by $a'^{\parallel a}$.
\end{definition}

An other characterization is the following:
\begin{lemma}
$a'$ is invertible along $a$ if and only if there exists $b\in S$ such that $ba'b=b$ and $b\cH a$, and in this case $a'^{\parallel a}=b$.
\end{lemma}

\begin{theorem}\label{thexist}
Let $a,a'\in S$. Then the following are equivalent:
\ben
\item $a'^{\parallel a}$ exists.
\item $a\leq_{\cR} aa'$ and $(aa')\grp$ exists.
\item $a\leq_{\cL} a'a$ and $(a'a)\grp$ exists.
\item $aa'a\cH a$ ($a'\in A(a)[\cH]$).
\item $a\leq_{\cH} aa'a$.
\een
In this case, $$b=a(a'a)^{\sharp}=(aa')^{\sharp}a.$$
\end{theorem}

In other words, $a'\in A(a)[\cH]$ if and only if $a'$ is invertible along $a$, or equivalently if and only if $a'$ is an associate of an element in the $\cH$-class of $a$. It was remarked in \cite{maryPatricio} that, even if Green's relation $\cH$ is not a congruence in general, invertibility modulo $\cH$ has a good behaviour with respect to $\cH$-classes. 

\begin{corollary}\label{corH}
$$aa'a\cH a\Longleftrightarrow \cH_a a' \cH_a =\cH_a$$ 
\end{corollary}

We will need the following characterization of reflexive inverses modulo $\cH$ (see \cite{mary}). Recall that $ab$ is a trace product if $ab\in \cR_a\cap \cL_b$, or equivalently, by a theorem of Clifford \cite{Clifford56}, if $\cL_a\cap \cR_b$ contains an idempotent.

\begin{theorem}\label{thtraceproduct}
The following statements are equivalent:
\ben
\item $a'\in V(a)[\cH]$;
\item $(a')^{\parallel a}$ exists and is an inverse of $a'$ ($(a')^{\parallel a}\in V(a')$);
\item $aa'$ and $a' a$ are trace product.
\een
\end{theorem}

In this case, it is known that the following equality hold \cite{Clifford56}: $$\cH_a\cH_{a'}=\cH_{aa'}=a\cH_{a'}=\cH_{a} a'.$$
Together with corollary \ref{corH} we get:
\begin{proposition}\label{prop-Hclasses}
If, $(a,a')$ is a regular pair modulo $\cH$, then $$\cH_a\cH_{a'}\cH_a=\cH_{a}=a\cH_{a'}a.$$
\end{proposition}


We now prove that all the classical properties of inverses remain true when working modulo $\cH$. Recall that for $a',a''\in A(a)$, then $a'a, aa' \in E(S)$ and $a'aa'' \in V(a)$. Also, if $a'\in V(a)$ then $a'a=aa' \iff a'\cH a$.

\begin{proposition}
Let $a',a''\in A(a)[\cH]$. Then $a'a, aa'\in E(S)[\cH]$ and $a'aa''\in V(a)[\cH]$. If moreover, $a'\in V(a)[\cH]$, then $a'a\cH aa' \iff a'\cH a$.
\end{proposition}

\proof
By theorem \ref{thexist}, $a'\in A(a)[\cH]$ implies $a'a$ and $aa'$ are group invertible. But $E(S)[\cH]=\cH(S)$ by Clifford's theorem. 
Second, let $a',a''\in A(a)[\cH]$, that is $aa'a\cH a$ and $aa''a\cH a$. Then by corollary \ref{corH}, $\cH_a a'\cH_a=\cH_a$. Applying this result to $a\in \cH_a$ and $aa''a\in \cH_a$ we get $aa'(aa''a)\in \cH_a$, and $a'aa''\in A(a)[\cH]$. To prove that $a'aa''\in V(a)[\cH]$, we have to prove that $a'aa''\leq_{\cH} a'aa''aa'aa''$. We prove only $\leq_{\cR}$, for the other one is symmetric. We use the inverses of $a'$ and $a''$ along $a$, $b'=(a')^{\parallel a}=a(a'a)\grp$ and $b''=(a'')^{\parallel a}=a(a'a)\grp$.
\begin{align*}
a'aa''&=a'\left(aa''\left(a''\right)^{\parallel a}\right)a''\\
&=a'aa''a(a''a)\grp a''\\
&=a'aa''\left(aa'\left(a'\right)^{\parallel a}\right)\left(a''a\right)\grp a''\\
&=a'aa''aa'a\left(a'a\right)\grp\left(a''a\right)\grp a''\\
&=a'aa''aa'\left(aa''\left(a''\right)^{\parallel a}\right)\left(a'a\right)\grp\left(a''a\right)\grp a''
\end{align*}
and $a'aa''\leq_{\cR}a'aa''aa'aa''$.\\
Finally, suppose $a'\in V(a)[\cH]$ and $a'a\cH aa'$. Pose $b=\left(a'\right)^{\parallel a}$. From $b\cR a$ we get $a'b\cR a'a$ ($\cR$ is a left congruence). From $aa'b=a$ we get $a\leq_{\cL} a'b \leq_{\cL}b\leq_{\cL}a$ and finally, since $a'a\cL a$,  $a'b\cH a'a$. Symmetrically, 
$ba'\cH aa'$, and finally $a'b\cH ba'$. But $a'b,ba'\in E(S)$. They are then idempotents in the same $\cH$-class, hence they are equal (corollary 1 in \cite{Clifford56}). It follows that $a'$ is the group inverse of $b$ ($(a,b)$ is a regular pair from theorem \ref{thtraceproduct}). Finally $\cH_b=\cH_a$ is a group  
that contains $a'$ (and $a\grp$). Conversely, if $a'\in V(a)[\cH]$ and $a'\cH a$ then $\cH_{aa'}=\cH_a\cH_{a'}=\cH_a.\cH_a=\cH_{a'}\cH_{a}=\cH_{a'a}$.
\endproof

Finally, we study classical types of semigroups with respect to our new definitions. First, we show that regularity is the same as regularity modulo $\cH$. 

\begin{lemma}\label{lemma_Hreg}
A semigroup is regular if and only $\forall a\in S, \;A(a)[\cH]\neq \emptyset$.
\end{lemma}

\proof
The implication is straightforward. For the converse, let $a'\in A(a)[\cH]$, that is $aa'a\cH a$. Then $(a')^{\parallel a}$ exists and $a=aa'(a')^{\parallel a}=aa'(aa')^{\sharp} a$ by theorem \ref{thexist}. Finally, $a$ is regular.
\endproof

Next theorem proves that inverse semigroups can also be defined in terms of inverses modulo $\cH$. Recall that $S$ is an inverse semigroup if elements admit a unique inverse $a^{-1}$ (or equivalently, $S$ is regular and $\forall a\in S,\; a', a''\in V(a)\Rightarrow a' = a''$). 

\begin{theorem}\label{thweakinverse}
let $S$ be a semigroup and $a\in S$ be a regular element. The following statements are equivalent:
\ben
\item $a', a''\in V(a)\Rightarrow a' = a''$;
\item $a', a''\in V(a)[\cH]\Rightarrow a'\cH a''$.
\een 
\end{theorem}

\proof
$\,$\\
\vspace{-.5cm}
\ben
\item[$(1)\Rightarrow (2)$] Let $a\in S,\; a', a''\in V(a)[\cH]$. Then by theorem \ref{thexist}, $a^{\parallel a'}$ and $a^{\parallel a''}$ are well defined, and by theorem \ref{thtraceproduct}, they are inverses of $a$. Since $S$ is inverse, they are equal. But by definition $a^{\parallel a'}\cH a'$ and $a^{\parallel a''}\cH a''$, and finally $a'\cH a''$.

\item[$(2)\Rightarrow (1)$] This is a result of Clifford (\cite{Clifford56}), any $\cH$-class contains at most one inverse. We prove it here for completeness. Let $a', a''\in V(a)$. Then they are $\cH$-inverses, and $a'\cH a''$. By cancellation properties, $a'\cL a''$ and $a'=a'aa'\Rightarrow  a''=a''aa'$, $a''\cR a'$ and $a''=a''aa''\Rightarrow  a'=a''aa'$. Finally $a'=a''aa'=a''$.

\een 
\vspace{-.5cm}
\endproof

As a direct corollary, we get that a semigroup $S$ is inverse if it satisfies $$\forall a\in S,\;V(a)[\cH]\neq \emptyset \hbox{ and }a', a''\in V(a)[\cH]\Rightarrow a'\cH a''.$$ 

One may wonder is this type of results is also true for other characterizations of inverse semigroups, or for orthodox semigroups. This is not the case as next section will show. 

\section{Completely inverse semigroups and $\cH$-orthodox semigroups}

It is well known that inverse semigroups can be alternatively defined as regular semigroups whose idempotents commute (this was actually both Vagner and Preston first definition, and the equivalence between the two notions was proved by Vagner \cite{Vagner52} and Liber \cite{Liber54}). Non-regular semigroups whose idempotents commute (resp. form a subsemigroup) have been studied for instance in \cite{Margolis87}, \cite{Fountain90}, \cite{Almeida92}, under various names, the latest being ``$E$-commutative semigroups'' (resp. ``$E$-semigroups''). Some of the results of the paper also apply to non-regular semigroups, and we try to avoid regularity assumption when it is not necessary.\\
Since we can consider group invertible elements as idempotents modulo $\cH$, we consider here semigroups whose group invertible elements $\cH$-commute (commutation modulo $\cH$ was introduced by Tully \cite{Tully73}). 
 
\begin{lemma}\label{lemma-CIisI}
Let $S$ be a semigroup such that $a,b\in \cH(S)\Rightarrow ab\cH ba$. Then $E(S)$ is commutative.
\end{lemma}  
 
\proof
Let $S$ be such a semigroup. Then group invertible elements hence idempotents $\cH$-commute. Cancellation properties shows that the products are actual. Indeed, $ef\cL fe$ and $eff=ef\Rightarrow fef=fe$. But also $ef\cR fe$ and $ffe=fe\Rightarrow fef=ef$  and idempotents commute. 
\endproof

Now, we show that the set of group invertible elements in an inverse semigroup is not necessarily $\cH$-commutative. Consider the inverse semigroup of partial injective transformations on $X=\{1,2\}$ $S=\mathcal{I}_X$. It is an inverse semigroup. However, $\cH(S)$ is not $\cH$-commutative. Indeed, define $\alpha$ by $\alpha(1)=1$ and $\beta$ by $\beta(1)=2,\; \beta(2)=1$. Then $\alpha$ and $\beta$ are group invertible ($\alpha\grp=\alpha$, $\beta\grp=\beta$), but $\beta\alpha$ and $\alpha\beta$ have not the same domain, hence are not $\cL$-related.  

These results show that $\cH(S)$ being a $\cH$-commutative set is a stronger notion that $E(S)$ being a commutative set. Since regular semigroups satsfying this property are necessarily inverse, and are characterized in terms of group invertible elements, we call them completely inverse semigroups (by analogy with completely regular semigroups). We keep the name of the non-regular analog for later on.

\begin{definition}
A semigroup $S$ is completely inverse if it is regular and $\cH(S)$ is a $\cH$-commutative set ($a,b\in \cH(S)\Rightarrow ab\cH ba$).
\end{definition}

Lemma \ref{lemma-CIisI} hence claims that a completely inverse semigroup is inverse.

In order to study precisely these semigroups, we need more properties on the set $\cH(S)$. This will be done through the study of semigroups whose group invertible elements form a subsemigroup and $\cH$-orthodox semigroups. Recall that an orthodox semigroup may be defined as a regular semigroup whose set of idempotents forms a subsemigroup. Lemma 1.3 in \cite{Reilly67} gives equivalent characterizations of orthodox semigroups.

\begin{lemma}[Lemma 1.3 in \cite{Reilly67}]
For a regular semigroup $S$ the following are equivalent:
\begin{enumerate}
\item $E(S)E(S)\subseteq E(S)$;
\item $e\in E(S)$, $(e,x)$ a regular pair implies $x\in E(S)$;
\item $(a,a')$ and $(b,b')$ regular pairs implies that $(ab,b'a')$ is a regular pair. 
\end{enumerate}
\end{lemma}


\begin{definition}
A semigroup $S$ is $\cH$-orthodox if it is regular and $a'\in V(a)[\cH], b'\in V(b)[\cH]\Rightarrow b'a'\in V(ab)[\cH]$.
\end{definition}

We investigate now if the other characterizations are valid.

\begin{theorem}\label{thorthodox}
Let $S$ be a semigroup. The following two conditions are equivalent.
\ben
\item $a'\in V(a)[\cH], b'\in V(b)[\cH]\Rightarrow b'a'\in V(ab)[\cH]$.
\item $\cH(S)$ is a semigroup.
\een
\end{theorem}

\proof
$\,$\\
\vspace{-.5cm}
\ben
\item[$1)\Rightarrow 2)$] Let $a,b\in \cH(S)$. Then $aa\grp\in V(a)[\cH]$ and $b\grp b\in V(b)[\cH]$. It follows that $(b\grp baa\grp)\in V(ab)[\cH]$, and in particular $(ab)^2=ab(b\grp b aa\grp)ab\cH ab$, hence $ab$ is group invertible by Clifford's theorem.

\item[$2)\Rightarrow 1)$]

Let $a,b\in S,\; a'\in V(a)[\cH],\; b'\in V(b)[\cH]$. We want to show that $b'a'\in V(ab)[\cH]$, that is $(ab).(b'a')$ and $(b'a').(ab)$ are trace product. Since $a,b$ and $a',b'$ play symmetric roles, it is sufficient to prove that $(ab).(b'a')\in \cR_{ab}\cap \cL_{b'a'}$ and using the opposite (dual) semigroup, it is sufficient to prove only the $\cL$ relation.
First, $bb'\cL b'\Rightarrow bb'a'\cL b'a'$. Pose $x=bb'a'a$. By theorem \ref{thexist}, $bb',a'a\in \cH(S)$ hence $x\in \cH(S)$. From $a'a\cR a'$ we get $bb'a'a\cR bb'a'$, and by cancellation $x=xx\grp x\Rightarrow bb'a'=xx\grp bb'a'$. Then
$bb'a'=x\grp x bb'a'=x\grp bb'a'a  bb'a'$. Finally, $b'a'\cL bb'a'\cL abb'a'$.



\een 
\vspace{-.5cm}
\endproof

The class of regular semigroups whose set of group invertible elements is a subsemigroup has indeed already been studied in a paper of Warne \cite{Warne76} under the name ``natural regular semigroups'', and appear also in some works of Yamada and Shoji \cite{Yamada79}, \cite{Yamada80}. The theorem claims that the class of $\cH-$orthodox semigroups, defined upon inverses modulo $\cH$, coincide with the class of natural regular semigroups. Non-regular semigroups whose set of group invertible elements is a subsemigroup appear in a paper of Birget, Margolis and Rhodes \cite{Birget90} under the name ``solid semigroups''. Remark that for such semigroups, $\cH(S)$ being a semigroup and a union of groups is completely regular, hence a semillatice of completely simple semigroups.\\

Surprisingly, these two conditions are not equivalent with the other possible, namely that $(h,x)$ is a regular pair modulo $\cH$ and $h\in \cH(S)$ implies $x\in \cH(S)$. For simplicity we will call a semigroup with this property a $\cH$-inverse-closed semigroup.

\begin{definition}
A semigroup $S$ is $\cH$-inverse-closed if $(h,x)$ is a regular pair modulo $\cH$, $h\in \cH(S)$ implies $x\in \cH(S)$.
\end{definition}

Such a semigroup has interesting properties. For instance:
\begin{lemma}
Let $S$ be a $\cH$-inverse-closed semigroup. 
\ben
\item $e,f\in E(S)\textit{ and } ef \textit{ regular}\Rightarrow ef\in \cH(S)$;
\item $a'\in A(a), b'\in A(b)\Rightarrow b'a'\in A(ab)[\cH]$;
\item $a'\in V(a), b'\in V(b)\Rightarrow b'a'\in V(ab)[\cH]$.
\een
\end{lemma}

\proof
$\,$\\
\vspace{-.5cm}
\ben
\item By a result of Howie and Lallement (lemma 1.1 in \cite{Howie66}), any regular product of idempotents $ef$ admits an inverse $p$ that is idempotent. Then $(p,ef)$ is a regular pair (\textit{a fortiori} modulo $\cH$) and $ef\in \cH(S)$. 
\item Let $a'\in A(a)$ and $b'\in A(b)$. Then $ab=aa'abb'b$, with $a'a,bb'\in E(S)$. By the first point, $h=a'abb'\in \cH(S)$. 
It follows that 
\begin{eqnarray*}
ab&=ahb=ah^2(h\grp)^3 h^2 b\\
&=aa'abb'a'abb'(h\grp)^3 a'abb'a'abb'b\\
&=(ab b'a'ab)(b'(h\grp)^3 a')(abb'a'ab).
\end{eqnarray*}
It follows that $ab\leq \cH ab b'a'ab$, that is $b'a'\in A(ab)[\cH]$.
\item Since being a regular pair is symmetric and $a$ and $b$ play symmetric roles, the result follows from $2)$.
\een
\vspace{-.5cm}
\endproof

\begin{theorem}
Let $S$ be a semigroup. Then:
\ben
\item $\cH(S)$ is a semigroup implies $S$ is $\cH$-inverse-closed;
\item $S$ is inverse implies $S$ is $\cH$-inverse-closed.
\een
\end{theorem}

\proof
$\,$\\
\vspace{-.5cm}
\ben
\item Assume $\cH(S)$ is a semigroup and let $h\in \cH(S), x\in V(h)[\cH]$. Then $xhx\cH x$. But $xhx=(xh)h\grp (hx)$ is the product if three group invertible elements by theorem \ref{thexist}, hence in the subsemigroup $\cH(S)$, and $x$ being in the $\cH$-class of a group invertible element is group invertible. 
\item Assume now that $S$ is inverse and let $h\in \cH(S), x\in V(h)[\cH]$. Then $(h,h\grp)$ and $(h,x)$ are regular pairs and by theorem \ref{thweakinverse}, $h\grp$ and $x$ are in the same $\cH$-class. By lemma \ref{lemma_Clifford}, $x$ is group invertible.
\een
\vspace{-.5cm}
\endproof

One may wonder why we do not call $\cH$-orthodox semigroups completely orthodox semigroups. The answer is that $\cH$-orthodox semigroups may not be orthodox. For instance, completely regular semigroups are obviously $\cH$-orthodox, but they are not all orthogroups (completely regular and orthodox). \\
On the other hand, note that orthodox or $\cH$-inverse-closed semigroups may not be $\cH$-orthodox. Consider the previous example of partial injective transformations on $X=\{1,2\}$ $S=\mathcal{I}_X$. It is inverse hence orthodox and $\cH$-inverse-closed. However, the domain and image of $\beta\alpha$ do not coincide, hence $\beta\alpha$ is not group invertible. This shows that $S=\mathcal{I}_X$ is not $\cH$-orthodox. \\

The remaining of this section is devoted to the study of the link between $\cH$-commutation and closedeness of the set of group invertible elements.
 
\begin{lemma}\label{lem-ROL}
Let $a,b\in \cH(S)$ such that $ab=bb\grp ab=abaa\grp$ and $ba=babb\grp =aa\grp ba$. Then $ab$ is group invertible and $(ab)\grp=b\grp a\grp$. 
\end{lemma}

\proof
We start by commutation:
\begin{eqnarray*}
(ab)b\grp a\grp &=bb\grp ab b\grp a\grp = b\grp (ba)b b\grp a\grp \\
&= b\grp (aa\grp ba) bb\grp a\grp = b\grp aa\grp (ba bb\grp) a\grp\\ 
&=b\grp aa\grp ba a\grp =b\grp a\grp (abaa\grp) = b\grp a\grp ab
\end{eqnarray*}
Note that by symmetry $baa\grp b\grp=a\grp b\grp ba$.
Second, we focus on inner invertibility.
\begin{equation*}
(ab)b\grp a\grp ab=b\grp a\grp abab=b\grp ba b =ab
\end{equation*}
We finally adress outer invertibility. 
\begin{equation*}
(ab)(b\grp a\grp) =(b\grp b ab)(b\grp a\grp)=b\grp (b ab b\grp) a\grp=bb\grp aa\grp
\end{equation*}
Symetrically, $aa\grp bb\grp=baa\grp b\grp =a\grp b\grp ba$.
Finally, 
\begin{equation*}
b\grp a\grp ab b\grp a\grp=b\grp baa\grp b\grp a\grp= b\grp b b\grp a\grp a a\grp=b\grp a\grp
\end{equation*} 
\endproof
 
\begin{theorem}\label{thROL}
Let $a,b\in \cH(S)$. If $ba\cH ab$, then $ab\in \cH(S)$ and $(ab)\grp=b\grp a\grp$.
\end{theorem}

\proof
Assume $ba\cH ab$, $a\cH a^2$ and $b\cH b^2$. Then by cancellation properties ($bb\grp ba=ba$,...)
$bb\grp ab=ab=abaa\grp$ and $aa\grp ba=ba=babb\grp$. We then apply the previous lemma.
\endproof

\begin{corollary}\label{corCIisHO}
Let $S$ be a semigroup. 
\ben
\item $\cH(S)$ is a $\cH$-commutative set implies $\cH(S)$ is a semigroup.
\item $S$ is completely inverse implies $S$ is $\cH$-orthodox.
\een
\end{corollary}

Note that in any inverse semigroup, the equality $(ab)^{-1}=b^{-1}a^{-1}$ gives that $b\grp a\grp$ is an inverse of $ab$ directly.
 
In order to study the converse implication, we need the following lemma.

\begin{lemma}
Let $a,b, ab, ba\in \cH(S)$ such that $(ab)\grp=b\grp a\grp$ and $(ba)\grp=a\grp b\grp$. Then $ab\cH ba$.
\end{lemma}

\proof
Assume $(ab)\grp=b\grp a\grp$ and $(ba)\grp=a\grp b\grp$. 
By cancellation properties,
\begin{equation*}
bb\grp b\grp a\grp=b\grp a\grp =b\grp a\grp aa\grp \Rightarrow bb\grp ab=ab=abaa\grp
\end{equation*}
But then 
\begin{equation*}
ab=bb\grp a\grp a ab=b (ab)\grp a bb\grp ab=bab(ab)\grp b\grp ab
\end{equation*}
and $ab\leq_{\cR} ba$. Conversely, 
\begin{equation*}
ab=ab bb\grp a\grp a=aba\grp a b(ab)\grp a=aba\grp (ab)\grp aba
\end{equation*}
Finally $ab\leq_{\cH} ba$. Since $a$ and $b$ play symmetric roles, also $ba\leq_{\cH} ab$ and finally, $ab\cH ba$.
\endproof

\begin{theorem}\label{thIHOisCI}
Let $S$ be a semigroup.
\ben
\item $\cH(S)$ is $\cH-$commutative if and only if $\cH(S)$ is a Clifford (inverse and completely regular) semigroup.
\item $S$ is completely inverse if and only if it is inverse and $\cH$-orthodox.
\een
\end{theorem}

\proof
The implication is lemma \ref{lemma-CIisI} and corollary \ref{corCIisHO}. For the converse, assume $\cH(S)$ is a Clifford semigroup and let $a,b\in \cH(S)$. Then $ab$ and $ba$ are group invertible by theorem \ref{thorthodox}. By unicity of the inverse in the inverse semigroup $\cH(S)$, $(ab)\grp=b\grp a\grp$ and $(ba)\grp=a\grp b\grp$. 
This implies that $ab\cH ba$ by the previous lemma, that is group invertible elements $\cH$-commute.
\endproof

Based on this theorem, we call semigroups whose group invertible elements $\cH$-commute $\cH$-Cliffordian semigroups.

\section{Commutation, Congruences and Crypticity}

Recall that a Clifford semigroup is either defined as an inverse and completely regular semigroup, or a regular semigroup with central idempotents (\cite{PetrichReilly-book}). Next theorem shows that this centrality is constitutive of completely inverse semigroups. 

First, we note that regular elements of the centralizer of $E(S)$ are always group invertible.
 
\begin{lemma}\label{lemma_inclusioncentral}
For any regular semigroup $S$, $Z(E(S))\subset \cH(S)$.
\end{lemma}

\proof
Let $a\in Z(E(S))$, $b\in V(a)$. Then $ab, ba\in E(S)$ hence $a=aba=a^2b=ba^2$ by commutation, and $a\leq_{\cH} a^2$. It follows that $a\in \cH(S)$.
\endproof

Now, we investigate the converse inclusion. For convenience, we restate the following result.

\begin{lemma}\label{lemma_incl}
Let $S$ be a semigroup such that $\cH(S)$ is a $\cH$-commutative set. Then $\cH(S)\subset Z(E(S))$.
\end{lemma}

\proof
Let $h\in \cH(S)$ and $e\in E(S)$. Then $eh\cH he$ and exists $x,y\in S^1,\, eh=xhe,\, he=ehy$. Computations give $eh=xhe=xhee=ehe=eehyehy=he$.
\endproof

\begin{corollary}\label{col-commut}
$S$ is a completely inverse semigroup if and only if $S$ is regular and $\cH(S)=Z(E(S))$.
\end{corollary}

\proof
From lemma \ref{lemma_inclusioncentral}, $Z(E(S))\subset \cH(S)$ and from lemma \ref{lemma_incl} $\cH(S)\subset Z(E(S))$ and the implication is proved.\\
For the converse, assume that $S$ is regular and $\cH(S)=Z(E(S))$. Then $S$ is inverse since $E(S)\subset \cH(S)$. But $Z(E(S))$ is always a semigroup, hence $S$ is $\cH$-orthodox. By theorem \ref{thIHOisCI}, $S$ is completely inverse.
\endproof

We now state a non-regular version of this result.
\begin{theorem}\label{th-commut}
$S$ is a $\cH$-Cliffordian if and only if $\cH(S)\subset Z(E(S))$.
\end{theorem}

\proof
The implication is lemma \ref{lemma_incl}.\\
For the converse, assume that $\cH(S)\subset Z(E(S))$. We first prove that $\cH(S)$ is a semigroup. Let $a,b\in \cH(S)$. Then $(ab)(aa\grp)=(aa\grp) ab=ab$ and symetrically $babb\grp=ba$. Also $ab=abb\grp b=b\grp b ab$ and $ba=a\grp a ba$. By lemma \ref{lem-ROL}, $b\grp a\grp $ is the group inverse of $ab$ and $\cH(S)$ is a semigroup. Finally $E(S)\subset \cH(S)\subset Z(E(S))$, and $\cH(S)$ is then a Clifford semigroup. We conclude by theorem \ref{thIHOisCI}.
\endproof

As a second corollary, we get that subsemigroups of completely inverse semigroup are completely inverse.

\begin{corollary}
Let $S$ be a completely inverse semigroup, and $T$ an inverse subsemigroup of $S$.
Then $T$ is completely inverse. 
\end{corollary}

\proof
First, $T$ is inverse (inverse semigroups form a variety) hence regular. Consider $Z_T(E(T))$. Since $T$ is regular $Z_T(E(T))\subset \cH_T(T)$. Conversely, $E(T)\subset E(S)$ and $Z(E(S))=\cH(S)$ implies that any element of $\cH(S)$ commutes with any idempotent in $E(T)$. Let $h\in \cH_T(T)$. Then $h\in \cH(S)$ and $h\in Z_T(E(T))$. Finally $T$ is regular and $Z_T(E(T))=\cH_T(T)$ and by corollary \ref{col-commut}, $T$ is completely inverse.
\endproof

To continue the study, we need the following definitions and theorems about $\cH$-commutative semigroups (\cite{Nagy-book}). 

\begin{definition}
A semigroup $S$ is $\cR$-commutative ($\cL$-commutative, $\cH$-commutative) if for all $a,b\in S$, exists $x\in S^1$, $ab=bax$ ($ab=xba$, $ab=bxa$). 
\end{definition}

At first sight, this definition of $\cH$-commutativity seems quite different from the one from Tully (\cite{Tully73} used in this paper: $\forall a,b\in S, ab\cH ba$. Next theorem (\cite{Nagy-book}) claims that they are equivalent for sets that are semigroups.

\begin{theorem}\label{th-Tully}
$\,$\\
\vspace{-.5cm}
\ben
\item A semigroup is $\cH$-commutative if and only if it is $\cR$ and $\cL$-commutative.
\item A semigroup $S$ is $\cR$-commutative ($\cL$-commutative, $\cH$-commutative) if and only if $\cR$ ($\cL$, $\cH$) is a commutative congruence (\textit{i.e.} it is a congruence and
the quotient semigroup is commutative).
\een
\vspace{-.5cm}
\end{theorem}

The interest of these definitions and theorem lies in the fact that the set $\cH(S)$ is a $\cH$-commutative semigroup when $S$ is completely inverse. Since the last theorem involves congruences, we start from some general facts. 

Congruences on semigroups have been studied intensively, for general semigroups and for classes of semigroups (regular, inverse, orthodox, completely regular, completely simple,...) as well. We recall the following results (see for instance \cite{Reilly67} and \cite{Scheiblich74}). An idempotent separating congruence is a congruence $\rho$ such that $(e,f)\in \rho$ and $e,f\in E(S)$ implies $e=f$. They can be characterized as those congruences such that $\rho \subset \cH$. From lemma 1.2 in \cite{Scheiblich74}, idempotent separating congruences on inverse semigroups are in $1:1$ correspondence with Kernels, inverse subsemigroups $K$ that are self-conjugate and satisfy $E(S)\subset K\subset Z(E(S))$. The correspondance is given by $(a,b)\in \rho_K\iff \{a^{-1}a=b^{-1}b\mbox{ and } ab^{-1}\in K\}$. Note that $E(S)$ and $Z(E(S))$ are Kernels when $S$ is inverse.

\begin{theorem}\label{thCIcong}
A semigroup $S$ is complelely inverse if and only if $S$ is inverse and $\cH$ is a congruence.
\end{theorem} 

\proof
$\,$\\
\vspace{-.5cm}
\begin{itemize}
\item[$\Rightarrow$]
Assume $S$ is completely inverse. Then $\cH(S)=Z(E(S))$ by corollary \ref{col-commut}, hence it is the Kernel of the congruence $\rho_{\cH(S)}$ with $(a,b)\in \rho_{\cH(S)}\iff \{a^{-1}a=b^{-1}b\mbox{ and } ab^{-1}\in \cH(S)\}$. We now prove that $\rho_{\cH(S)}=\cH$. First, $\rho_{\cH(S)}\subset \cH$ since $\rho_{\cH(S)}$ is idempotent separating (lemma 1.2 in \cite{Scheiblich74}). Let $a\cH b$. First, $a\cL b$ is equivalent to $a^{-1}a=b^{-1}b$ in any inverse semigroup (lemma \ref{lemma_L}). Also, from $\cH_b\cH_{b^{-1}}=\cH_{bb^{-1}}$ (proposition \ref{prop-Hclasses}) we get  $ab^{-1}\cH bb^{-1}$ hence $ab^{-1}\in \cH(S)$.
Finally $a\cH b\iff (a,b)\in \rho_{\cH(S)}$ and $\cH$ is a congruence on $S$.
\item[$\Leftarrow$] Conversely, assume $S$ is inverse and $\cH$ is a congruence. The homomorphic image of an inverse semigroup is always inverse (\cite{Preston54}), hence $S/\cH$ is inverse. But idempotents in $S/\cH$ are the $\cH$-classes of elements of $\cH(S)$. It follows that idempotents in $S/\cH$ commute if and only if $\cH(S)$ is $\cH$-commutative, and the result follows.
\end{itemize}
\vspace{-.5cm}
\endproof

Semigroups such that $\cH$ is a congruence are sometimes called $\cH$-compatible semigroups \cite{Yamada79}, or following Reilly and Petrich cryptic semigroups.
The previous theorem shows that among cryptic semigroups, the inverse ones admit a rich structure (for instance, the set of group invertible elements is a Clifford semigroup).

Actually, we can improve a little bit the only if condition using a result of Lallement \cite{Lallement67} for regular semigroups: ``if $\cH$ is a congruence on the semigroup $\cH(S)$, then it is a congruence on $S$''. 

\begin{corollary}
A semigroup $S$ is complelely inverse if and only if $S$ is inverse and 
$\forall a,b,c\in \cH(S),\; a\cH b\Rightarrow ca\cH cb \text{ and } ac\cH bc$.
\end{corollary}

Actually, this statement holds for non-regular semigroups.
  
\begin{theorem}\label{th-nonregcong}
Let $S$ be a semigroup.
\ben
\item If $\forall a,b,c\in \cH(S),\; a\cH b\Rightarrow ca\cH cb \text{ and } ac\cH bc$, and $S$ is an $E$-semigroup then $\cH(S)$ is a semigroup ($S$ is solid);
\item $S$ is $\cH$-Cliffordian ($\cH(S)$ is a $\cH$-commutative set) if and only if $S$ is an $E$-commutative semigroup and 
$\forall a,b,c\in \cH(S),\; a\cH b\Rightarrow ca\cH cb \text{ and } ac\cH bc$.
\een
\end{theorem}

\proof
$\,$\\
\vspace{-.5cm}
\ben
\item Let $a,b\in \cH(S)$ and pose $e=aa\grp,\, f=bb\grp$. Then $a\cH e$, $b\cH f$. By the congruence property, it follows that $ab\cH eb$ and $eb\cH ef$, and finally $ab\cH ef$. But $ef$ is idempotent since $S$ is an $E$-semigroup hence $\cH(S)$ is a semigroup. 
\item If $S$ is $\cH$-Cliffordian, then by lemma $E(S)$ is commutative. Also $\cH(S)$ is a semigroup by corollary \ref{corCIisHO}, and by theorem \ref{th-Tully}, $\cH$ is a congruence on $\cH(S)$. \\
For the converse, we use the previous construction. For $a,b\in \cH(S)$ we pose $e=aa\grp,\, f=bb\grp$. Then $ab\cH ef$ and $ba\cH fe$. But $ef=fe$ hence $ab\cH ba$, and $\cH(S)$ is a $\cH$-commutative semigroup.
\een
\vspace{-0.5cm}
\endproof

We deduce directly from theorem \ref{thCIcong} the following corollary.

\begin{corollary}\label{cor-quotient}
Let $S$ be a semigroup. Then $S$ is completely inverse if and only if it is cryptic and $S/\cH$ is inverse. In this case, $S/\cH$ is combinatorial (hence completely inverse).
\end{corollary}

\proof
If $S$ is completely inverse, then $\cH$ is a congruence by the previous theorem and $S/\cH$, as the homomophic image of an inverse semigroup, is inverse. Conversely, if $\cH$ is a congruence and idempotents in $S/\cH$ commute, then $\cH(S)$ is $\cH$-commutative. Also any element of $S$ is regular modulo $\cH$ hence regular, and the semigroup is completely inverse. \\
Let $H_a,H_b\in S/\cH$. Then by lemma \ref{lemma-L} $$H_a\cH_{S/\cH} H_b\iff H_a H_{a^{-1}}=H_b H_{b^{-1}} \text{ and } H_{a^{-1}} H_a = H_{b^{-1}} H_b.$$ But $H_x H_{x^{-1}}=H_{xx^{-1}}$ for all $x\in S$ hence $aa^{-1}\cH bb^{-1}$ and since they are idempotent, they are equal. It follows that $a=aa^{-1} a=bb^{-1}a$ and $b=bb^{-1}b=aa^{-1}b$, hence $a\cR b$. The same arguments in the opposite semigroup give $a\cL b$, and finally $\cH_a=\cH_b$.
\endproof

\begin{corollary}\label{cor-selfc}
Let $S$ be a completely inverse semigroup and $(a,a')$ be a regular pair modulo $\cH$, that is $a'\in V(a)[\cH]$. Then $\cH_{a'}\cH(S)\cH_{a}\subset \cH(S)$.
\end{corollary}

\proof
By corollary \ref{cor-quotient}, $S/\cH$ is inverse and $E(S/\cH)$ is self-conjugate. 
\endproof

Once again, the regularity assumption is not necessary, but the proof is more involved.

\begin{lemma}\label{lemma-selfc}
Let $S$ be a semigroup such that $\cH(S)$ is $\cH$-commutative, and let $(a,a')$ be a regular pair modulo $\cH$. Then $a'\cH(S)a\subset \cH(S)$.
\end{lemma}

\proof
By theorem \ref{corCIisHO}, $\cH(S)$ is a $\cH$-commutative semigroup.
Let $a'\in V(a)[\cH]$ and $h\in \cH(S)$. Then $haa'\cH aa'h$ since $\cH(S)$ is a $\cH$-commutative semigroup and $aa'\in \cH(S)$. Also, $\cH$ is a congruence on $\cH(S)$ hence $h\cH h^2 \Rightarrow haa'\cH h^2aa'\cH haa'h\cH aa'h$. Also $h(aa')\grp\cH h(aa')\grp h\cH (aa')\grp h$. Then exists $x,y\in S^1$, $h(aa')\grp=xh(aa')\grp h$ and $haa'=yhaa'h$. it follows that $$a'ha=a'haa'(a')^{\parallel a}=a'haa'(aa')\grp a=a'h(aa')\grp aa'a.$$ Using $x$ then $y$ we get 
$$a'ha=a'xh(aa')\grp h aa'a=a'xh(aa')\grp yh aa'ha.$$
But also $haa'\cH haa'haa'$ since $\cH(S)$ is a semigroup, and exists $z\in S^1$, $haa'=z haa'haa'$. Finally 
$$a'ha=(a'xh(aa')\grp y)z haa'haa'ha$$ and $a'ha\leq_{\cL}(a'ha)(a'ha)$. Symmetrically, $a'ha\leq_{\cR}(a'ha)(a'ha)$ and finally $a'ha\in \cH(S)$.
\endproof

Finally, recall that in any regular inverse semigroup, the maximal congruence included in $\cH$ is of the form 
$(a,b)\in \mu \iff \{\forall e\in E(S), a^{-1}ea \cH b^{-1}eb\}$ or equivalently $(a,b)\in \mu \iff  \{a^{-1}a = b^{-1}b \mbox{ and } ab^{-1}\in Z(E(S))\}$ \cite{Howie64}. Define a new equivalence relation by $(a,b)\in \nu \iff \{\forall h\in \cH(S), a^{-1}ha \cH b^{-1}hb\}$.

\begin{corollary}
On a completely inverse semigroup $S$, $\rho_{\cH(S)}=\nu=\mu=\cH$.
\end{corollary}

\proof
First, $\mu\subset \cH$ since it separates idempotents. But $\cH$ is a idempotent separating congruence on $S$ completely inverse hence $\cH\subset \mu$, and we get $\cH=\mu$. By corollary \ref{col-commut} $\cH(S)=Z(E(S))$ hence by the previous description of $\mu$, $\rho_{\cH(S)}=\mu$. Also, $E(S)\subset \cH(S)\Rightarrow \nu\subset \mu$. Let $(a,b)\in \mu$ and $h\in \cH(S)$. Then $hh\grp$ is idempotent and $h\cH hh\grp$. But relation $\cH$ is a congruence on the completely inverse semigroup $S$, hence
$a^{-1}ha\cH a^{-1}hh\grp a \cH b^{-1}hh\grp b\cH b^{-1}hb$ and $(a,b)\in \nu$. \endproof


\section{Conclusion}
It is worth summarizing all the previous characterizations of completely inverse semigroups. A semigroup $S$ is completely inverse if it verifies one of the following equivalent conditions:

\ben
\item $S$ is regular and $\cH(S)$ is $\cH$-commutative.
\item $S$ is regular and $\cH(S)$ is a Clifford semigroup.
\item $S$ is regular and $\cH(S)=Z(E(S))$.
\item $S$ is inverse and $\cH$-orthodox.
\item $S$ is inverse and $\cH$ is a congruence ($S$ is cryptic inverse).
\item $\cH$ is a congruence and $S/\cH$ is inverse. 
\een

Also, we have the following implications between the different types of semigroups studied in this paper.
\begin{figure}[h!]
$$\xymatrix{
& \cH-\text{commutative and regular}  \ar@{=>}[d] \ar@{=>}[ld]& Cryptic\\
\text{Completely Regular} \ar@{=>}[d]  & \text{Completely Inverse}\ar@{=>}[rd] \ar@{=>}[ru] \ar@{=>}[ld] & + \ar@{=>}[l]  \\
 \cH-\text{Orthodox} \ar@{=>}[rd] & + \ar@{=>}[u]&\text{Inverse}\ar@{=>}[ld] \ar@{=>}[d]  \\  
& \cH-\text{Inverse-closed} & \text{Orthodox} 
}\label{figimplications}$$
\end{figure}

We finally give a few example of completely inverse semigroups, $\cH$-Cliffordian semigroups and $\cH$-orthodox semigroups. 

\ben
\item Clifford semigroups (in particular regular commutative semigroups) are completely inverse.
\item Completely regular semigroups are $\cH$-orthodox.
\item by theorem \ref{th-nonregcong}, any cryptic orthodox semigroup is $\cH$-orthodox.
\item The bicylic semigroup is completely inverse. Here $\cH(S)=E(S)$ since $a\cH b\iff a=b$.
\item More generally, any combinatorial ($a\cH b\Rightarrow a=b$) inverse semigroup is completely inverse since $\cH(S)=E(S)$.
\item A fundamental inverse semigroup is completely inverse if and only if it is combinatorial.
\item Consider the inverse semigroup of partial injective transformations on $X=\{1,2,3,4\}$, $\mathcal{I}_X$. Let $S$ be its regular subsemigroup with the following elements: $S=\{0,e,f,a,b,g,h\}$ with $0$ the empty function, $e:1\mapsto 1$, $f:2\to 2$ , $g:(3,4)\mapsto (3,4)$, $a:1\mapsto 2$, $b:2\mapsto 1$ and $h:(3,4)\mapsto (4,3)$. $S$ is a completely inverse semigroup but neither completely regular nor combinatorial. Indeed, $E(S)=\{0,e,f,g\}$ and $\cH(S)=\{0,e,f,g,h\}$. 
\item Let $G$ be a group, $g\in G$ fixed and pose $S=G\cup\{a\}$, $a^2=g^2$, $ah=gh \forall h\in G$,  $ha=hg \forall h\in G$. Then $S$ is $\cH$-Cliffordian (here $\cH(S)=G=\cH_g$ is a group). It is not completely inverse since $a$ is not a regular element.
\een 

As a final comment, note that the class of completely inverse semigroups is not a variety of $(2,1)$-algebras (algebras with the two operations of multiplication and inversion). Indeed, we have proved that a subsemigroup of a completely inverse semigroup is completely inverse. However, this is not true for the homomorphic image. As a counterexample, take $X$ an inverse not completely inverse semigroup, and consider $(S=\cF_X,f)$ the free inverse semigroup on $X$. Then by the universal property of the free inverse semigroup, for $i:X\to X$ the identity map there exists a (unique) homomorphism $h:\cF_X\to X$ such that $hf=i$, and $X$ is the homomorphic image of $\cF_X$. But $\cF_X$ is combinatorial hence completely inverse (Reilly, Lemma 1.3 in \cite{Reilly72}), whereas $X$ is not.

\bibliographystyle{alpha}

\end{document}